\documentclass[11pt]{amsart}
\usepackage{graphicx}
\usepackage[usenames]{color}
\usepackage{amssymb,amscd}
\newcommand\Z{\mathbb Z}
\newcommand\Q{\mathbb Q}
\newcommand\R{\mathbb R}
\newcommand\A{{\mathbb H}^2}
\newcommand\Rinf{R_{\infty}}
\newcommand\G{\Gamma_n}
\newcommand\Gs{\Gamma(S)}
\newcommand\T{\Pi_{i=1}^k T_i}
\newcommand\s{(n_1,n_2, \cdots n_k)}
\newcommand\sprime{(m_1,m_2, \cdots ,m_l)}

\newcommand\pt{\Pi_{i=1}^k T_i}

\newcommand\ptp{\Pi_{i=1}^k T'_i}
\newcommand\ph{\varphi}

\newtheorem{theorem}{Theorem}[section]

\newtheorem{lemma}[theorem]{Lemma}
\newtheorem{corollary}[theorem]{Corollary}
\newtheorem{definition}[theorem]{Definition}

\title{Twisted conjugacy and quasi-isometry invariance for generalized solvable Baumslag-Solitar groups}
\author{Jennifer Taback and Peter Wong}
\address{Department of Mathematics,
Bowdoin College, Brunswick, ME 04011} \email{jtaback@bowdoin.edu}
\thanks{The first author acknowledges support from
NSF grant DMS-0437481, and would like to thank Kevin Whyte for many useful conversations about this paper. 
We thank the anonymous referee for helpful comments.}
\address{Department of Mathematics, Bates College, Lewiston, ME 04240} \email{pwong@bates.edu}

\begin{document}

\maketitle

\begin{abstract}
We say that a group has property $\Rinf$ if any group automorphism has an infinite number of twisted conjugacy classes.  Fel'shtyn and Gon\c calves prove that the solvable Baumslag-Solitar groups $BS(1,m)$ have property $\Rinf$. We define a solvable generalization $\Gamma(S)$ of these groups which we show to have property $\Rinf$.  We then show that property $\Rinf$ is geometric for these groups, that is, any group quasi-isometric to $\Gamma(S)$ has property $\Rinf$ as well.
\end{abstract}

\section{Introduction}
The celebrated Lefschetz Fixed Point Theorem asserts that a self map $f$ of a compact manifold $M$ has a fixed point when the Lefschetz number $L(f)$ does not vanish.  Additionally, under these conditions, any map homotopic to $f$ possesses a fixed point as well.  Since the Lefschetz number is merely an algebraic count of the number of fixed points, the converse of this theorem does not hold in general.

A more subtle homotopy invariant is the Nielsen number $N(f)$, which is defined to be the number of essential fixed point classes of $f$.  Often the Nielsen number is used to prove the converse of the Lefschetz Fixed Point Theorem.  Since the Nielsen number $N(f)$ is a lower bound for the number of fixed points of any element in the homotopy class of $f$, the ability to compute $N(f)$ may yield the existence of multiple fixed points.

When $f$ is a selfmap of a compact manifold of dimension at least three, the Nielsen number $N(f)$ is a sharp lower bound on the number of fixed points of any element in the homotopy class of $f$.  That is,  there exists a map $f'$ homotopic to $f$ with exactly $N(f)$ fixed points. Nielsen fixed point theory has proven useful in nonlinear analysis, differential equations, and dynamical systems. Thus the computation of $N(f)$ is a central issue in Nielsen fixed point theory.

The calculation of the Nielsen number is in general quite difficult.  The first computational result was due to W. Franz (1943) who studied selfmaps of a lens space.  He showed that all fixed point classes of such maps have the same fixed point index, so that 
\begin{enumerate}
\item $L(f)=0$ implies $N(f)=0$, in which case $f$ is deformable to a fixed point free map, and 
\item $L(f)\ne 0$ implies $N(f)=\#Coker (1-f_{*_1})$ where $f_{*_1}$ is the induced homomorphism on the first integral homology. 
\end{enumerate}
When the Lefschetz number $L(f)$ is nonzero, the cardinality of $Coker (1-f_{*_1})$ is exactly equal to the Reidemeister number $R(f)$.  The Reidemeister number is defined to be the cardinality of the set of $f_{\sharp}$-twisted conjugacy classes in the fundamental group, where $f_{\sharp}$ is the map induced by $f$ on the fundamental group of $M$.  Thus in certain cases, the computation of $N(f)$ reduces to the computation of $R(f)$, which is an easier quantity to calculate. 

In 1963, B. Jiang identified a large class of spaces, now known as Jiang spaces, satisfying the so-called Jiang condition for which either (1) or (2) holds for all selfmaps.  A space $X$ is said to satisfy the {\it Jiang condition} if the evaluation subgroup $Im(ev_*:\pi_1(X^X,1_X) \to \pi_1(X,x_0))$ coincides with the fundamental group $\pi_1(X,x_0)$ where $ev_*$ is induced by the evaluation map at the point $x_0\in X$. Recently, Jiang-type results have been proven for certain non-Jiang spaces. In particular, conditions were given so that $N(f)=R(f)$, {\it provided} $R(f)$ is finite. In general, when $R(f)$ is finite it provides an upper bound for $N(f)$, even when $N(f)\ne R(f)$.

We now ask the following question.  For which families of (fundamental) groups $\pi$ is the Reidemeister number infinite for all $\varphi: \pi \rightarrow \pi$, or at least for all automorphisms $\varphi$, excluding the equality of the Reidemeister number and the Nielsen number?  Analogously, which families of (fundamental) groups do not have $R(\varphi) = \infty$ for all $\varphi$, allowing the possible computation of the Nielsen number by equality with the Reidemeister number for certain selfmaps $\varphi$?

In 1994, A. Fel'shtyn and R. Hill  \cite{fh} related the Reidemeister torsion with the dynamical zeta function of an endomorphism $\varphi$ whose coefficients are the Reidemeister numbers of the iterates $\varphi^n$ of $\varphi$. They conjectured that every injective endomorphism $\varphi:\pi \to \pi$ of a finitely generated torsion-free group $\pi$ of exponential growth must have infinite number of twisted conjugacy classes.  In \cite{ll}, G. Levitt and M. Lustig implicitly showed that $R(\varphi)$ is infinite for all automorphisms $\varphi$ of a finitely generated torsion-free non-elementary Gromov hyperbolic group $\pi$. Fel'shtyn \cite{f} extended \cite{ll} by removing the torsion-free hypothesis and pointed out, using the work of Z. Sela \cite{s} on co-Hopficity, that the Fel'shtyn-Hill conjecture holds for finitely generated torsion-free non-elementary and freely indecomposable Gromov hyperbolic groups.  In 2003, D. Gon\c calves and P. Wong \cite{gw1} constructed a class of polycyclic groups of exponential growth which are mapping tori of Anosov automorphisms that are not Gromov hyperbolic for which automorphisms with finite Reidemeister number exist. It is natural to look for classes of groups for which every automorphism has infinite Reidemeister number.

A finitely generated group $\pi$ is said to have property $R_{\infty}$ if $R(\varphi)=\infty$ for all automorphisms $\varphi$ of $\pi$. With this definition, torsion-free non-elementary Gromov hyperbolic groups have property $R_{\infty}$ \cite{f}. Recently, A. Fel'shtyn and D. Gon\c calves \cite{fg} showed that the class of Baumslag-Solitar groups $BS(m,n)$ have property $R_{\infty}$ except in the case $m=1=n$.  On the other hand, D. Gon\c calves and P. Wong \cite{gw2} recently classified all such finitely generated abelian groups $G$ for which the (restricted) wreath product $G \wr \mathbb Z$ has property $R_{\infty}$.  They show that the lamplighter groups $\Z_n \wr \Z$ have property $\Rinf$ only when $2$ or $3$ divides $n$.

In this paper, we extend the class of groups known to have property $R_{\infty}$ using a variety of techniques. We first define a solvable generalization of $BS(1,n)$ which arises in an analogous geometric context, and a further generalization which does not.  We give two proofs that these groups have property $\Rinf$, one algebraic and one which relies on the geometry of the groups.  We then show that property $\Rinf$ is geometric for these groups, meaning that any group quasi-isometric to $\Gamma(S)$ has property $\Rinf$ as well.

\section{Background on twisted conjugacy classes}

Let $\ph:\pi \to \pi$ be a group endomorphism and  let $\pi$ act on itself via $\sigma \cdot \alpha\mapsto \sigma \alpha \ph(\sigma)^{-1}$ for $\sigma, \alpha \in \pi$. The orbits of this action are the {\it Reidemeister classes} of $\ph$ or the $\ph$-twisted conjugacy classes. Denote by $R(\ph)$ the cardinality of the set $\mathcal R(\ph)$ of $\ph$-twisted conjugacy classes.  This number $R(\ph)$ is called the {\em Reidemeister number} of $\ph$. When $\ph$ is the identity, $R(\ph)$ is the number of conjugacy classes of elements in $\pi$.

Consider the following commutative diagram where the rows are short exact sequences of groups and the vertical arrows are group homomorphisms, that is, $\ph|A = \ph'$ and ${\overline \ph}$ is the quotient map induced by $\ph$ on $C$:

\begin{equation}\label{short-exact}
\begin{CD}
    1   @>>> A    @>{i}>>  B @>{p}>>      C @>>> 1 \\
    @.  @V{\ph'}VV      @V{\ph}VV   @V{\overline \ph}VV @.\\
    1   @>>> A    @>{i}>>  B @>{p}>>      C @>>> 1 
 \end{CD}
\end{equation}

This induces a short exact sequence of sets and corresponding functions ${\hat i}$ and ${\hat p}$:
$$
\mathcal R(\ph') \stackrel{\hat i}{\to} \mathcal R(\ph) \stackrel{\hat p}{\to} \mathcal R(\overline \ph) \to 1
$$
where if $\bar 1$ is the identity element in $C$, we have $\hat{i}(\mathcal R(\ph'))={\hat p}^{-1}([\bar 1])$,  and $\hat p$ is onto. Here, the functions $\hat i$ and $\hat p$ are the canonical functions induced by $i$ and $p$ respectively. The fixed subgroup $Fix \overline \ph=\{c\in C|\overline \ph(c)=c\}$ acts on $\mathcal R(\ph')$ via $c\cdot [\alpha] \mapsto [b\alpha \ph(b)^{-1}]$ where $p(b)=c$. 
It is easily seen that $[b\alpha \ph(b)^{-1}] \in {\mathcal R}(\ph')$, since $$p(b \alpha \ph(b)^{-1}) = p(b) p(\alpha) p(\ph(b)^{-1}) =  c {\overline \ph}(c)^{-1} = cc^{-1} = \bar 1.$$
Moreover, $\hat i([\alpha_1])=\hat i([\alpha_2])$ iff $[\alpha_1]$ and $[\alpha_2]$ belong to the same orbit of the action of $Fix \overline \ph$.

The following result is straightforward and follows from more general results discussed in  Theorem 1 of \cite{w}.

\begin{lemma}\label{lemma:reid}
Given the commutative diagram labelled \eqref{short-exact} above, 
\begin{enumerate}
\item if $R(\overline \ph)=\infty$ then $R(\ph)=\infty$, and 

\item if $C$ is finite and $R(\ph')=\infty$ then $R(\ph)=\infty$.
\end{enumerate}
\end{lemma}

For certain endomorphisms $\ph$ of group extensions, the Reidemeister number $R(\ph)$ is determined by the Reidemeister numbers of the restriction and the induced homomorphism, namely $R(\ph')$ and $R({\overline \ph})$.

\begin{lemma}\label{lemma:reid2}
Consider the commutative diagram \eqref{short-exact} and let $\tau_{\alpha} \beta=\alpha^{-1}\beta \alpha$. Suppose for any $\bar \alpha \in C$, we have $Fix (\tau_{\bar \alpha}\bar \ph) =\{\bar 1\}$. If $\ph_0:B\to B$ induces the same commutative diagram as \eqref{short-exact}, that is, $\ph_0 i=i \ph'$ and $\bar \ph p=p \ph_0$, then $R(\ph)=R(\ph_0)$.
\end{lemma}
\begin{proof}
It follows from Corollary 1 of \cite{w} that if $R(\ph)<\infty$ then $R(\ph)=\sum_{[\bar \alpha]} R(\tau_{\alpha}\ph')$ where $p(\alpha)=\bar \alpha$ and $[\bar \alpha]\in \mathcal R(\bar \ph)$. Since $\sum_{[\bar \alpha]} R(\tau_{\alpha}\ph')$ only depends on $\ph'$ and $\bar \ph$, it follows that $R(\ph)=R(\ph_0)$. 

Similarly, if $R(\ph)=\infty$ but $R(\ph_0)< \infty$, then the argument above shows that $R(\ph_0)$ depends only on $\ph'$ and $\bar \ph$.  Since $\ph'$ and $\bar \ph$ are identical whether we consider $\ph$ or $\ph_0$, we would have to conclude that $R(\ph)<\infty$.  Thus we must have $R(\ph_0)= \infty$ as well.
\end{proof}

The following will be used in the subsequent sections and applied in particular to semidirect products.

\begin{corollary}\label{cor:reid}
Consider the commutative diagram \eqref{short-exact} with the following additional assumptions:
\begin{enumerate}
\item the extension splits,
\item $C=\mathbb Z^k=\langle t_1,..., t_k|t_it_j=t_jt_i\rangle$ for $k\ge 1$, and
\item the endomorphisms are automorphisms.
\end{enumerate}
Let $\ph_0:B\to B$ be an automorphism which induces the same commutative diagram as \eqref{short-exact} and has the property that $\ph_0(t_i)\in \mathbb Z^k$.  Then $R(\ph)=R(\ph_0)$.
\end{corollary}
\begin{proof}
The existence of such an automorphism $\varphi_0$ follows from the fact that the extension splits. If $R(\bar \ph)=\infty$ then the statement follows from (1) of Lemma \ref{lemma:reid} because $\varphi$ and $\varphi_0$ induce the same automorphism $\bar \varphi$. Suppose that $R(\bar \ph)<\infty$. Since $C$ is abelian $\tau_{\bar \alpha}\bar \ph=\bar \ph$ for any $\bar \alpha \in C$. It is well known that $R(\bar \ph)<\infty$ iff $Fix \bar \ph=\{\bar 1\}$. The result then follows from Lemma \ref{lemma:reid2}.
\end{proof} 

\section{Background on quasi-isometries and quasi-actions}

A quasi-isometry is a map between metric spaces which distorts distance by a uniformly bounded amount, defined precisely as follows.
\begin{definition}
Let $X$ and $Y$ be metric spaces.  A map $f: X \rightarrow Y$ is a $(K,C)$-quasi-isometry, for $K \geq 1$ and $C \geq 0$ if
\begin{enumerate}
\item $\frac{1}{K} d_X(x_1,x_2) - C \leq d_Y(f(x_1),f(x_2)) \leq K d_X(x_1,x_2) +C$ for all $x_1,x_2 \in X$.
\item For some constant $C'$, we have $Nhbd_{C'}(f(X)) = Y$.
\end{enumerate}
\end{definition}

There is a notion of {\em coarse inverse} for a quasi-isometry; quasi-isometries $f$ and $g$ are coarse inverses if both compositions $f \circ g$ and $g \circ f$ are a bounded distance from the identity.  If $X$ and $Y$ are path connected metric spaces, there is a standard procedure called ''connecting the dots" which allows one to alter a quasi-isometry by a uniformly bounded amount, replacing it with a continuous map satisfying the quasi-isometry conditions above.  This procedure is described, for example, in \cite{sw}.  We will assume that our quasi-isometries below have been altered in this way.

Milnor, Svarc and Effremovich \cite{m} observed independently that if a finitely generated group acts properly discontinuously and cocompactly by isometries on a proper geodesic metric space $X$, then $G$ is quasi-isometric to $X$.  In particular, if $G$ is finitely generated and finitely presented, the Cayley complexes of a group with respect to two different generating sets are quasi-isometric.  With this notion, we can associate, up to quasi-isometry, a {\em geometric model} to a finitely generated and finitely presented group, namely the Cayley complex of the group which is quasi-isometric to it.  We always use the word metric to measure distance in the group.  

We say that two groups are {\em abstractly commensurable} if they contain isomorphic finite index subgroups.  Since any finite index subgroup of a finitely generated group is quasi-isometric to the entire group, abstractly commensurable groups are quasi-isometric.

The set of all self quasi-isometries of a finitely generated group $G$ is denoted $QIMap(G)$.  We form equivalence classes consisting of all quasi-isometries of $G$ which differ by a uniformly bounded amount.  This set of equivalence classes is called the {\em quasi-isometry group} of $G$ and denoted $QI(G)$.  

If $G$ is a finitely generated group and $X$ is a proper geodesic metric space, we define a {\em quasi-action} of $G$ on $X$ to be a map $\Psi: G \rightarrow QIMap(X)$ satisfying the following properties for some constants $K \geq 1$ and $C \geq 0$.
\begin{enumerate}
\item For each $g \in G$, the element $\Psi(g)$ is a $(K,C)$-quasi-isometry of $X$.
\item $\Psi(Id)$ is a bounded distance from the identity, that is, $\Psi(Id)$ is the identity in $QI(X)$.
\item $\Psi(g) \Psi(h)$ is a bounded distance from $\Psi(gh)$, for all $g,h \in G$.
\end{enumerate}
It is clear that a quasi-action induces a homomorphism from $G$ into $QI(X)$.  It is important to note that all quasi-isometries $\Psi(g)$ have the same quasi-isometry constants. 

If $f: G \rightarrow X$ is a quasi-isometry, then the quasi-action $\Psi: G \rightarrow QIMap(X)$ is defined by $\Psi(g) = f \circ L_g \circ f^{-1}$, where $L_g$ is left multiplication by $g \in G$.  In this situation, we can characterize the kernel of this quasi-action, that is, those $g \in G$ so that $\Psi(g)$ is a uniformly bounded distance from the identity quasi-isometry, as the virtual center of the group $G$.  The {\em virtual center $V(G)$} is the subgroup of $G$ consisting of all elements whose centralizers have finite index in $G$, that is, $V(G)=\{g\in G| [G:C(g)]<\infty\}$.  This characterization of the kernel of the quasi-action follows from the following lemma, which is straightforward.  A proof is given in \cite{two}, and relies on the fact that $g \in G$ moves all points of $G$ a uniformly bounded distance under left multiplication if and only if $g$ has finitely many conjugates.

\begin{lemma}\label{virtual-center} 
Let $G$ be a finitely generated group with a quasi-action on a Cayley complex $X$.  Then the virtual center of $G$ consists of those elements that move all points of $X$ a uniformly bounded distance $B$, that is, $ \ d(gx,x) \leq B$ for all $x \in X$.
\end{lemma}

\begin{corollary}\label{cor:kernel}
Let $G$ be a finitely generated group quasi-isometric to a Cayley complex $X$.  Let $\Psi: G \rightarrow QIMap(X)$ be a quasi-action.  Then $Ker(\Psi) = V(G)$.
\end{corollary}

\section{The solvable Baumslag-Solitar groups and a generalization}
The solvable Baumslag-Solitar groups $BS(1,m)$ have presentation $\left< a,t | tat^{-1} = a^m \right>$, and can be expressed as the group extension
$$1 \rightarrow \Z[ {\hbox{$\frac{1}{m}$}}] \rightarrow BS(1,m) \rightarrow \Z \rightarrow 1.$$
The large-scale geometry of these groups was studied by Farb and Mosher, who describe the geometric model for these groups as the warped product of a tree and the real line  \cite{fm1}.  The metric on this complex is such that when a height function is put on the tree, and $l$ is a line in the tree on which the height function is strictly increasing, $l \times \R$ is quasi-isometric to a hyperbolic plane.  For more details on this complex, we refer the reader to \cite{fm1}.

We consider a solvable generalization of $BS(1,m)$, which we define in two stages.  The initial generalization arises in a geometric context, but easily extends  to a strictly algebraic generalization.  When we view the group $PSL_2(\Z[\frac{1}{p}])$, for $p$ a prime, as a subgroup of $PSL_2(\R) \times PSL_2(\Q_p)$ under the diagonal embedding, we get a natural action of $PSL_2(\Z[\frac{1}{p}])$ on $\A \times T$ where $T$ is the Bruhat-Tits tree associated to $PSL_2(\Q_p)$.  Under this action, the stabilizer of a point at infinity is exactly $BS(1,p^2)$, which is commensurable with $BS(1,p)$.

We can analogously consider the group $PSL_2(\Z[\frac{1}{n}])$, where $n$ is any positive integer, which acts on the product of $\A$ with the product of the Bruhat-Tits trees for $PSL_2(\Q_{p_i})$.  Here, $p_i$ are the prime factors of $n$.   The stabilizer of a point at infinity under this action is the upper triangular subgroup of $PSL_2(\Z[\frac{1}{n}])$, which we denote $\G$.  This group is also defined by the group extension $1 \rightarrow \Z[\frac{1}{n}] \rightarrow \G \rightarrow \Z^k \rightarrow 1$.  We use the following presentation of $\G$:
$$\G \cong \left< a,t_1,t_2, \cdots ,t_k | t_it_j = t_jt_i, \ i \neq j,\  t_i^{-1}at_i = a^{p_i} \right>.$$
Note that when $k=1$ this is exactly $BS(1,p)$.  The geometric model for these groups is topologically the product of $\R$ with the product of trees on which the group acts, and is detailed in \cite{tw}.  

The groups $\G$ generalize immediately to a larger class of groups $\Gs$ which do not arise in the same geometric context, but have similar algebraic properties to $\G$.  Let $S = \{n_1,n_2, \cdots ,n_k\}$ be a set of pairwise relatively prime positive integers.  We define $\Gs$ to be the group
$$\Gs \cong \left< a,t_1,t_2, \cdots ,t_k | t_it_j = t_jt_i, \ i \neq j,\  t_i^{-1}at_i = a^{n_i} \right>$$
and we use this presentation throughout the paper.
We will use the following definition of $\Gs$ as a group extension in the proofs below:
$$1 \rightarrow A \rightarrow \Gs \stackrel{\pi} \rightarrow \Z^k \rightarrow 1$$
where the map $\pi: \Gs \rightarrow \Z^k$ is a ``component-wise height function"  where the $i$-th coordinate in $\pi(g)$ is the sum of the exponents of the generator $t_i$ in a word representing $g \in \Gs$.  As with $BS(1,m)$ and $\G$, the kernel $A$ of $\pi$ is the normal closure of the subgroup generated by the element $a$ in the presentation above.  These groups have geometric models analogous to those for $\G$; see \cite{tw} for more details.

\section{The groups $\Gs$ have property $\Rinf$}
\label{RinfGamman}

We first prove that the solvable groups $\Gs$ which are generalizations of $BS(1,n)$ have property $\Rinf$, that is, any group automorphism has infinitely many twisted conjugacy classes.  We include two proofs of this fact, one algebraic, and one which relies on the geometry of the groups in a fundamental way.

We begin by expressing $\Gs$ as the group extension $1 \rightarrow A \rightarrow \Gs \stackrel{\pi} \rightarrow \Z^k \rightarrow 1$, and showing that the kernel $A$ is invariant under an automorphism of $\Gs$.  This allows us to obtain commutative diagrams from the short exact sequence describing $\Gs$.  We then use either the presentation or the geometry of the group to prove that the induced quotient map on $\Z^k$ is factor preserving, and indeed the identity on each factor.  The theorem then follows from Lemma \ref{lemma:reid}.

Since the Baumslag-Solitar groups $BS(1,m)$ represent a special case of $\Gs$, the algebraic proof given below uses techniques similar to those of \cite{fg}.  The groups $\Gs$ are also semidirect products $A \rtimes \Z^k$, as are the solvable Baumslag-Solitar groups when $k=1$.  In the geometric proof that $\Gs$ has property $\Rinf$ we utilize the geometry of the group via the fact that any group automorphism is also a quasi-isometry.  We use two facts about quasi-isometries of $\Gs$ from \cite{tw}, stated below for completeness.  First, we need that any quasi-isometry preserves the two factors of the
geometric model, namely $\R$ and $\Pi_{i=1}^{k} T_i$.  Second, the fact that the quasi-isometry group can be described explicitly.  These are both extensions of the initial results of Farb and Mosher \cite{fm1} concerning quasi-isometries of the solvable Baumslag-Solitar groups.

\begin{theorem}[\cite{tw} Corollary 4.3]
\label{thm:tw}
Consider the two groups $\Gamma_1 = \Gamma\s$
and $\Gamma_2 = \Gamma\sprime$ where $(n_i,n_j) = (m_i,m_j) = 1$ for $i \neq j$, and a quasi-isometry $f: \Gamma_1
\rightarrow \Gamma_2$ between them.  Then:
\begin{enumerate}
\item
$k=l$.
\item
$f$ induces a quasi-isometry $f_T: \pt \rightarrow \ptp$.
\end{enumerate}
\end{theorem}

\begin{theorem}[\cite{tw}, Theorem 1.5] \label{thm:tw3}
Let $S = \{n_1,n_2, \cdots ,n_k\}$ be a set of pairwise disjoint integers.  The quasi-isometry group of $\Gs$ is isomorphic to the product
$$Bilip (\R) \times Bilip (\Q_{n_1}) \times \cdots \times Bilip(\Q_{n_k}).$$
\end{theorem}

When $\Gamma(S)$ is written as a group extension $1 \rightarrow A \rightarrow \Gs \rightarrow \Z^k \rightarrow 1$, the subgroup $A$ is the normal closure of the generator $a$.  To see that this subgroup is characteristic, we note that it can be viewed as the normal closure of the set of elements which are conjugate to a (nontrivial) power of themselves.  The next lemma follows immediately from this characterization.

\begin{lemma}\label{lemma:char}
Write $\Gs$ as
$$1 \rightarrow A \rightarrow \Gs \rightarrow \Z^k \rightarrow 1$$
and let $\varphi \in Aut(\Gs)$.  Then $A$ is invariant under $\varphi$.
\end{lemma}

Using Lemma \ref{lemma:char} we obtain commutative diagrams from the short exact sequence defining $\Gs$, the first step towards proving that $\Gs$ has property $\Rinf$.

\begin{theorem}\label{thm:Rinfinity}
The group $\Gs$ has property $\Rinf$.
\end{theorem}

\begin{proof}
We see from Lemma \ref{lemma:char} that $A$ is invariant under $\varphi \in Aut(\Gs)$, and obtain the following commutative diagram.

\begin{equation}\label{general-exact}
\begin{CD}
    1 @>>> A   @>>>  \Gs @>>>    \Z^k @>>> 1 \\
    @.     @V{\varphi'}VV  @V{\varphi}VV   @V{\overline \varphi}VV @.\\
    1 @>>> A        @>>>  \Gs @>>>    \Z^k @>>> 1 
 \end{CD}
\end{equation}
Since $\Gs$ is a semidirect product, there is a map from $\Z^k$ back to $\Gs$, and thus can view $\Z^k$ as generated by $t_1, \cdots t_k$ in the presentation above for $\Gs$.  

Since $\Gs$ is a semidirect product, using Corollary \ref{cor:reid} we replace $\varphi$ by $\varphi_0 \in Aut(\Gs)$ with the following properties:
\begin{enumerate}
\item $\varphi_0' = \varphi'$
\item $\varphi_0(t_i)$ is a word in $t_1, \cdots t_k$, and
\item $R(\varphi) = R(\varphi_0)$.
\end{enumerate}
Using property (2) above, write $\overline{\varphi}_0(t_i) = w_i(t_1, \cdots t_k)$ as a word in the generators $t_1, \cdots ,t_k$.  Without loss of generality, we write $\varphi$ for $\varphi_0$.

We now finish the proof using the geometry of the group.  An algebraic proof is also given below.

Let $X$ denote the geometric model of $\Gs$. Then $X$ is topologically the product of $\R$ with $k$ trees, $\T$. It follows from Theorem \ref{thm:tw} (\cite{tw}, Corollary 4.3) that any quasi-isometry of $X$ preserves the factors $\R$ and $\T$. If $f$ is a quasi-isometry of $X$, we write $F|_{\R}$ and $f|_{T}$ for $f$ restricted to $\R$ and $\T$ respectively.
 The description of the quasi-isometry group $QI(\Gamma(S))$ given in Theorem \ref{thm:tw3} above implies that a quasi-isometry $f$ of $\Gamma(S)$ is a uniformly bounded distance from a quasi-isometry $\psi = (\psi_{\R},\psi_1,\psi_2, \cdots , \psi_k)$ of $\Gamma(S)$, where $\psi_{\R}: \R \rightarrow \R \in Bilip(\R)$, and $\psi_i \in Bilip(\Q_{n_i})$ induces a quasi-isometry on $T_i$ which we also denote $\psi_i$.  Since any automorphism of $\Gamma(S)$ is a quasi-isometry of $X$ as well, $\varphi \in Aut(\Gamma(S))$ is a uniformly bounded distance from such a quasi-isometry $\psi = (\psi_{\R},\psi_1,\psi_2, \cdots , \psi_k)$.
%
%

We now describe a particular $\Z^k$ subcomplex of $X$, namely the restriction of the map $\Gs \rightarrow X$ to the subgroup generated by the elements $\{ t_i\}$.  We will denote this particular $\Z^k$ by $\Z_*^k$.

Since $\varphi$ is an automorphism of $\Gs$ and induces vertical maps on the short exact sequence describing $\Gs$, in particular on $\Z^k \cong \left< t_1, \cdots ,t_k \right>$, and $\overline{\varphi}(t_i) = w_i(t_1, \cdots t_k)$, we see that $\varphi$ acting on the complex $X$ preserves $\Z_*^k$, i.e. $\varphi(\Z_*^k) = \Z_*^k$.  We now consider $\psi(\Z_*^k)$.

\begin{lemma} Let $\psi \in QI(\Gamma(S))$ and $\Z_*^k$ be defined as above.  Then
$\psi(\Z_*^k)  = \Z_*^k$.
\end{lemma}

\begin{proof}
As proved in Theorem 8.1 of \cite{fm1}, combining the maps $(\psi_{\R},\psi_i)$ yields a quasi-isometry of an $\R \times T_i$ subspace of $X$, which is metrically a Baumslag-Solitar subcomplex corresponding to $BS(1,n_i)$.  When a height function is imposed on each tree, and we take a line $l \subset T_i$ on which the height function is strictly increasing, the product $\R \times l$ is metrically a hyperbolic plane.

Let $l_i$ be the line in $\Z_*^k \subset X$ which is the image of $\{t_i^n\}$ for all $n \in \Z$.  Then the height function on $T_i$ is strictly increasing on $l_i$ by construction.  This line is contained in a Baumslag-Solitar subcomplex of the form $\R \times T_i$ of $X$, and in fact is contained in a unique hyperbolic plane in this subcomplex.  It follows from Proposition 4.1 of \cite{fm1} that quasi-isometries of $BS(1,m)$ coarsely preserve the hyperbolic planes in the complex which are of the form $\R \times l$, where $l$ is a line on which the height function is strictly increasing.  Thus, the image of this hyperbolic plane under $\psi$ is within a uniformly bounded distance of a unique hyperbolic plane in $\R \times T_i$ which is also of the form $\R \times l_i'$, where $l_i'$ is a line on which the height function is strictly increasing. When projected to the tree factor, this image plane yields a line $l_i'$ in $T_i$ which is the image of the original line $l_i$.  Since $\psi_T$ is a product of quasi-isometries, we conclude that $\psi(\Z_*^k)$ is the $\Z^k$ subspace of $X$ with axes $\{l_1', \cdots ,l_k'\}$ which intersects $(\psi_{\R},\psi_i)(Id)$ for each $i$.  

Since $d_X(\varphi(x),\psi(x)) \leq B$ for all $x \in X$, we must have $l_i' = l_i$.  If this is not the case, we see that there are points on $l_i$ and $l_i'$ which get arbitrarily far apart, since these are distinct lines in a tree.
\end{proof}

We have shown that $\varphi(\Z_*^k) = \psi(\Z_*^k) = \Z_*^k$. We now show that $\varphi$ is factor preserving as a map of $\Z_*^k$. 

We use coordinates $(a_1, \cdots ,a_k)$ on $\Z_*^k$ to represent $\Sigma_{i=1}^k a_it_i$.  Since $\psi_T$ is factor preserving, and $\psi$ preserves $\Z_*^k$, the map $\psi|_{\Z_*^k}$ is given by a diagonal matrix $D$ with diagonal entries $\{d_i\}$.  We know that $\varphi|_{\Z_*^k}$ is given by some matrix $M$ with integer entries and determinant $\pm 1$.

Let $v = (v_1,v_2, \cdots ,v_k) \in \Z^k$ be an eigenvector for $M$, with eigenvalue $\lambda$.  Then $D(Nv) = (d_1Nv_1, \cdots d_kNv_k)$ and $M(Nv) = \lambda Nv$ for any $N \in \Z$. Since $d_{\Z^k}(D(Nv),M(Nv)) \leq B$ for all $N \in \Z$, we see that for large values of $N$ we must have $d_i = \lambda$ for all $i = 1,2, \cdots ,k$.  Moreover, the product of the eigenvalues of $M$ is $\pm 1$, so each $d_i = \pm 1$, that is, $D$ is $\pm Id$.  

If $D$ is the identity matrix, we consider any vector $w \in \Z^k$.  Using the fact that $D$ and $M$ are a uniformly bounded distance apart, we write $d_{\Z^k}(Nw,M(Nw)) = d_{\Z^k}(Nw,NM(w)) = Nd_{\Z^k}(w,M(w))$.  This quantity cannot be uniformly bounded independent of $w$ unless we have $d_{\Z^k}(w,M(w)) = 0$ for all $w$, i.e. $M$ is also the identity matrix.  The same argument shows that if $-D$ is the identity matrix, so is $-M$.  

However, if $\varphi|_{\Z^k} = -Id$, so that $\varphi(t_i) = -t_i$, it is easy to see that the relator $t_iat_i^{-1} = a_i^{m_1}$ is not preserved under $\varphi$.  Thus $\varphi|_{\Z^k}$ is the identity map, and the theorem follows from Lemma \ref{lemma:reid}.
\end{proof}

We can alternately provide an algebraic proof of Theorem \ref{thm:Rinfinity}, as follows.

{\em Alternate proof of Theorem \ref{thm:Rinfinity}.} Follow the proof given above, replacing $\varphi$ by $\varphi_0 \in Aut(\G)$ with the following properties:
\begin{enumerate}
\item $\varphi_0' = \varphi'$
\item $\varphi_0(t_i)$ is a word in $t_1, \cdots t_k$, and
\item $R(\varphi) = R(\varphi_0)$.
\end{enumerate}
Using property (2) above, write $\overline{\varphi}_0(t_i) = w_i(t_1, \cdots t_k)$.  Without loss of generality, we write $\varphi$ for $\varphi_0$.  We show that indeed $\varphi(t_i)$ must be $t_i$, and thus the map induced on $\Z^k$ is the identity.  The theorem then follows from Lemma \ref{lemma:reid}.

We know that $\varphi(t_i)$ is a word in the generators $\{ t_j \}$, which we can write as $\Pi_{l=1}^k t_i^{\epsilon_i}$ since the $t_i$ commute with each another.  Since $t_i$ conjugates $a$ to $a^{m_i}$, we know that $\varphi(t_i)$ conjugates $\varphi(a)$ to $\varphi(a)^{m_i}$ as well.  Since the kernel of $\pi$ is characteristic, and so left invariant under any group automorphism, we must have $\phi(a) \in Ker(\pi)$.  Thus the sum of the exponents of $t_i$ in $\varphi(a)$ must be zero for all $i$.  Using the relators of $\Gs$, we can write $\varphi(a) = (\Pi_{l=1}^k t_i^{-e_i})a^f(\Pi_{l=1}^k t_i^{e_i})$.  It is easy to see, using the group relators $t_iat_i^{-1} = a^{m_i}$ that unless $\varphi(t_i) = t_i$ that $\varphi(t_iat_i^{-1})$ does not reduce to $\varphi(a)^{m_i}$.
\qed

\section{Property $\Rinf$ is a quasi-isometry invariant for $\Gs$}
\label{Rinfqi}

We now show that any group quasi-isometric to $\Gs$ also has property $\Rinf$.  In the proofs below, $S = \{ n_1,n_2, \cdots ,n_k \}$ is always a set of pairwise relatively prime integers.  This includes the special case of groups quasi-isometric to $BS(1,n)$.  Groups quasi-isometric to $\Gs$ are explicitly described in the following theorem of Taback and Whyte \cite{tw}.

\begin{theorem}[\cite{tw}, Theorem 1.6] \label{thm:tw1}
Let $\Gamma'$ be any finitely generated group quasi-isometric to $\Gs$.  Then there exist integers $m_1,m_2, \cdots ,m_k$ with each $m_i$ a rational power of $n_i$ as well as a finite normal subgroup $F$  of $\Gamma'$ so that $\Gamma'/F$ is isomorphic to a cocompact lattice in $Iso(X) \cong \R \rtimes (Sim(\Q_{m_1}) \times \cdots \times Sim(\Q_{m_k}))$.
\end{theorem}

In the special case when all $n_i$ are prime, and $\Gs$ is the upper triangular subgroup of the form $\G$, Theorem \ref{thm:tw1} can be refined to include a more specific description of the quotient group.

\begin{theorem}[\cite{tw}, Theorem 1.3] \label{thm:tw2}
Let $\Gamma'$ be any finitely generated group quasi-isometric to $\Gamma_n$.  Then there is a finite normal subgroup $F$ of $\Gamma'$ so that $\Gamma'/F$ is abstractly commensurable to $\Gamma_n$, meaning that $\Gamma'/F$ and $\Gamma_n$ have isomorphic subgroups of finite index.
\end{theorem}

If we specialize further to $\G = BS(1,m)$, Theorem \ref{thm:tw2} reduces to the rigidity theorem obtained by Farb and Mosher in \cite{fm2}.

The proof that any group quasi-isometric to $\Gs$ also has property $\Rinf$ does not explicitly use Theorems \ref{thm:tw1} or \ref{thm:tw2}.  Rather, we use the fact that such a group has a quasi-action on a product of trees, which can be quasi-conjugated, using the following theorem of \cite{msw}, to an isometric action on a product of possibly different trees.

\begin{theorem}[\cite{msw}, Theorem 1]
\label{thm:msw}
If $G \times T \rightarrow T$ is a cobounded quasi-action of a group $G$ on a bounded valence bushy tree $T$, then there is a bounded valence, bushy tree $T'$, an isometric action $G \times T' \rightarrow T'$ and a quasi-conjugacy $f: T \rightarrow T'$ from the action of $G$ on $T'$ to the quasi-action of $G$ on $T$.
\end{theorem}

To prove that a finitely generated group $G$ quasi-isometric to $\Gs$ has property $\Rinf$, we will need the fact that $G$ has finite center. We prove this using the quasi-action of $G$ on $X=\R \times \pt$ induced by the quasi-isometry $f: G \rightarrow \Gamma(S)$, and the characterization of the kernel of this quasi-action as the virtual center $V(G)$ given in Corollary \ref{cor:kernel}.  Since the center $Z(G)$ is contained in the virtual center $V(G)$, we need only to show that $V(G)$ is finite.  

To understand the kernel of this quasi-action, we make the following definition, and quote a lemma of Whyte \cite{wh}.  A space $X$ is called {\em tame} if for every $K \geq 1$ and $C> 0$ there is a constant $R>0$ so that every $(K,C)$-quasi-isometry is either an infinite distance from the identity or moves no point more than $R$.  

\begin{lemma}[\cite{wh}]\label{lemma:tame}
If $X$ is a cocompact, non elementary Gromov-hyperbolic space, then $X$ is tame.
\end{lemma}

To prove that $G$ has finite center we use the fact that each tree $T_i$ in the product $\pt$ is a cocompact, non elementary Gromov-hyperbolic space, and hence tame.

\begin{lemma}\label{lemma:finitecenter}
Let $S = \{ n_1,n_2, \cdots ,n_k\}$ be a set of pairwise relatively prime integers and $G$ any finitely generated group quasi-isometric to $\Gs$.  Then the center of $G$ is finite.
\end{lemma}

\begin{proof}
Let $\Theta: G \rightarrow QIMap(X)$ be the quasi-action of $G$ on $X$ which arises because $G$ is quasi-isometric to $\Gamma(S)$, and hence to $X$.  It follows from Theorem \ref{thm:tw3} that each $\Theta(g)$ is a uniformly bounded distance from a product map $\psi = (\psi_{\R},\psi_1,\psi_2, \cdots \psi_k)$, where $\psi_{\R} \in Bilip(\R)$, and $\psi_i \in Bilip(\Q_{n_i})$ induces a quasi-isometry of $T_i$ which we also denote $\psi_i$.  Each $T_i$ satisfies the conditions of Lemma \ref{lemma:tame} and thus is tame.

Let $g \in Ker(\Theta)$.  Then $\Theta(g)$ is a uniformly bounded distance from the identity quasi-isometry of $X$.  If $g \in Ker(\Theta)$ then $\Theta(g)$ is not a bounded distance from a quasi-isometry of $X$ which permutes the factors of $\pt$.  
All quasi-isometries $\Theta(g)$ have the same quasi-isometry constants $K$ and $C$ by definition, and we may assume that all quasi-isometries $\psi_i$ have the same quasi-isometry constants $K$ and $C$ as well.  

Since $\Theta(g)$ is defined using left multiplication $L_g$ on $G$, there are constants $R'$ and $R''$ depending on the tameness constant $R$, the uniformly bounded distance between $\Theta(g)$ and $\psi$, and the quasi-isometry constants $K$ and $C$ so that if $d_G(h,L_g(h)) \leq R'$ for all $h \in G$, then $d_X(x, \Theta(g) \cdot x) \leq R''$ for all $x \in X$, and hence when we restrict to $T_i$, points in this tree will be moved at most $R$ by $\psi_i$.  

Since $G$ is finitely generated, there are a finite number of $g \in G$ which move all points at most $R'$. These $g$ will lie in $Ker(\Theta)$ by construction.  For $g \in G$ which move the identity further than this constant $R'$, the induced quasi-isometry of $T_i$ will move points further than $R$, and so these quasi-isometries will be at infinite distance from the identity, so $g \notin Ker(\Theta)$.  Thus $Ker(\Theta)$ is finite, so the virtual center $V(G)$ and hence the center $Z(G)$ are also finite.
\end{proof}

We now prove the following theorem.

\begin{theorem}\label{thm:Rinfinity-qi}
Let $S = \{ n_1,n_2, \cdots ,n_k\}$ is a set of pairwise relatively prime integers and $G$ any finitely generated group quasi-isometric to $\Gs$. Then $G$ has property $\Rinf$.
\end{theorem}

\begin{proof}
Let $G$ be quasi-isometric to $\Gamma(S)$, equivalently to $X=\R \times \pt$.  Then there is a quasi-action $\Theta: G \rightarrow QIMap(X)$.  Identifying quasi-isometries which differ by a uniformly bounded amount, we view $\Theta(g)$ as an element of the quasi-isometry group $QI(X)$.  According to Theorem \ref{thm:tw}, any quasi-isometry $\Theta(g)$ of $X$ induces a quasi-isometry of $\pt$ which we denote $\Theta(g)|_T$.  For some $g \in G$, the quasi-isometry $\Theta(g)|_T$ may permute the tree factors of $\pt$. 

There is a finite group $F$ of permutations of these tree factors, and we obtain a surjective homomorphism $\pi: QI(G) \rightarrow F$.  If $g \in  Ker(\pi)$, then 
$\Theta(g)|_T$ acts as a quasi-isometry of $\pt$ which does not permute the tree factors.  Since the automorphism group of $G$ is a subgroup of the quasi-isometry group, by restriction, we may consider $\pi: Aut(G) \rightarrow F$, with kernel $K$.

We now show that $\varphi \in Aut(G)$ induces vertical maps on the short exact sequence 
$$1 \rightarrow K \rightarrow Aut(G) \stackrel{\pi} \rightarrow F \rightarrow 1.$$
Consider the homomorphism $i:G \rightarrow Aut(G)$ which maps $g \in G$ to the automorphism $\theta_g(\gamma) = g\gamma g^{-1}$ for $\gamma \in G$.  Note that $Ker(i)=Z(G)$, the center of $G$. For any automorphism $\varphi \in Aut(G)$, we have a commutative diagram of short exact sequences
\begin{equation}\label{center-exact}
\begin{CD}
    1     @>>>    Z(G)  @>>>       G        @>i>>  i(G)   @>>>     1 \\
    @.            @VV{\varphi'}V             @VV{\varphi}V  @VV{\overline {\varphi}}V \\
    1     @>>>    Z(G)  @>>>       G        @>i>>  i(G)   @>>>     1 
 \end{CD}
\end{equation}

For any $\ph \in Aut(G)$, there is an induced homomorphism $\eta_{\ph}:Aut(G) \to Aut(G)$ given by $\eta_{\ph}(\psi)=\ph \psi \ph^{-1}$. Note that the following diagram is commutative for any automorphism $\ph$  and $\eta_{\ph}|_{i(G)}=\overline {\ph}$.
\begin{equation*}
\begin{CD}
    G        @>i>>  Aut(G) \\
    @V{\varphi}VV  @VV{\eta_\varphi}V \\
    G        @>>i>  Aut(G) 
 \end{CD}
\end{equation*}
From this diagram, we see that if we can show that $R(\eta_{\varphi}) = \infty$, it follows that $R( \varphi) = \infty$.  We first show that when we consider $\pi: Aut(G) \rightarrow F$, the kernel $K = Ker(\pi)$ is invariant under $\eta_{\ph}$ for any automorphism of $\ph$ of $G$.

\begin{lemma}
\label{lemma:Kinvariant}
Let $\pi: Aut(G) \rightarrow F$ have kernel $K$.  Then $K$ is invariant under $\eta_{\ph}$ for any $\varphi \in Aut(G)$.
\end{lemma}

\begin{proof}
If $\varphi, \ \psi \in Aut(G)$, then $\eta_\varphi(\psi) = \varphi \psi \varphi^{-1} $. We must show that for any $\psi \in K = Ker(\pi)$, that $\eta_\varphi(\psi) \in K$, that is, $\pi(\eta_\varphi(\psi)) = 1 \in F$. But $\pi(\eta_\varphi(\psi)) = \pi (\varphi \psi \varphi^{-1}) = \pi(\ph)\pi(\psi)\pi(\ph^{-1}) = 1$ since $\psi \in K$ and $\pi$ is a homomorphism.
\end{proof}

We now use Lemma \ref{lemma:Kinvariant} to show that $K \cap i(G)$ is invariant under $\pi$ as well.  Recall that $i(G)= \{ \theta_g:G \rightarrow G | \theta_g(\gamma) = g \gamma g^{-1}  \text{ for } \gamma \in G \}$.  We must show that if $\theta_g \in K \cap i(G)$, then $\eta_\varphi(\theta_g) \in K\cap i(G)$ for any $\varphi \in Aut(G)$.  We know that $\eta_\varphi(\theta_g) \in K$ from Lemma \ref{lemma:Kinvariant}.

For any $\alpha \in G$, we have that 
$$\eta_\varphi(\theta_g)(\alpha) = \varphi(\theta_g \varphi^{-1}(\alpha)) = \varphi(g \varphi^{-1}(\alpha) g^{-1}) = \varphi(g) \alpha \varphi(g^{-1}) = \theta_{\varphi(g)}.$$
Thus $\eta_\varphi(\theta_g) \in i(G)$, and $K\cap i(G)$ is invariant under $\eta_{\ph}$.

We have now shown that we can induce vertical maps on the following short exact sequence.
\begin{equation}\label{general-exact}  
\begin{CD}
    1 @>>> K\cap i(G)   @>>>  i(G) @>\pi>>    F @>>> 1 \\
    @.     @V{\eta'_{\ph}}VV  @V{\eta_{\ph}}VV   @VVV @.\\
    1 @>>> K\cap i(G)        @>>>  i(G) @>\pi>>    F @>>> 1 
 \end{CD}
\end{equation}

We finish the proof by showing that any $\eta'_{\ph}: K\cap i(G) \rightarrow K\cap i(G)$ satisfies $R(\eta'_{\ph}) = \infty$, and the theorem then follows from Lemma \ref{lemma:reid}.

\begin{lemma}
\label{lemma:phi'}
Let $\pi: i(G) \rightarrow F$ be as above, $\varphi \in Aut(G)$ and $\eta'_{\ph}: K\cap i(G) \rightarrow K\cap i(G)$ the map induced on $K\cap i(G)$ as above.  Then $R(\eta'_{\ph}) = \infty$.  
\end{lemma}

\begin{proof}
Since $Z(G)$ is finite and $i(G) \cong G/Z(G)$, it follows from the corollary in Section 4.24 of \cite{dlH} that $G$ and $i(G)$ are quasi-isometric.
Thus there is a quasi-action of $i(G)$ on $\R \times \pt$. Restricting this quasi-action to $K\cap i(G)$, we have $\Theta':K\cap i(G) \rightarrow QIMap(\R \times \pt)$.  For every $g \in K\cap i(G)$, the quasi-isometry $\Theta'(g)$ of $\R \times \pt$ induces a quasi-isometry $\Theta'(g)|_T$ of $\pt$ by Theorem \ref{thm:tw}.  Moreover, since we have restricted from $i(G)$ to $K\cap i(G)$, we know that the quasi-isometry $\Theta'(g)|_T$ induces a map on $\pt$ which does not permute the tree factors.  We also know from Theorem \ref{thm:tw3} that $\Theta'(g)|_T$ is a bounded distance from a quasi-isometry $\psi$ of $\pt$ which is factor preserving, that is $\psi = (\psi _1,\psi _2, \cdots ,\psi _k)$ where each $\psi _i \in Bilip(\Q_{n_i})$ and hence induces a quasi-isometry $\psi_i$ of $T_i$.  Applying Theorem 1 of \cite{msw}, stated above as Theorem \ref{thm:msw}, the action of each $\psi _i$ on $T_i$ can be quasi-conjugated to give an {\em isometric} action on a possibly different tree $T_i'$.

Thus we get an isometric action of  $K\cap i(G)$ on the product of the new trees, $\ptp$.  Each $g \in K\cap i(G)$ yields a quasi-isometry  $\Theta'(g)|_T$ of $\ptp$ which is a uniformly bounded distance from the product map $\psi=(\psi_1,\psi_2, \cdots ,\psi_k)$ where each $\psi_i$ now has an isometric action on $T_i'$.  

Since each quasi-isometry $\Theta'(g)$ for $g \in K\cap i(G)$ preserves the point at infinity, when we quasi-conjugate to obtain the isometric action of $\psi_i$ on $T_i'$, the point at infinity is still preserved. Thus we can fix a height function on each $T_i'$, that is, a map $h_i:T_i' \rightarrow \Z$ as follows.  Fix a basepoint  vertex at height zero, and the height of any other vertex is defined relative to this fixed vertex.

Consider the height in $T_i'$ of the image of the basepoint of the tree $T_i'$ under the quasi-isometry $\psi_i:T_i' \rightarrow T_i'$.  We use this to determine a map 
$\xi_i:K\cap i(G) \rightarrow \Z$, where $\xi_i(g)$ is the height of the image of a vertex at height zero in $T_i$ under $\Theta'(g)|_{T_i'}$.  Combining these maps $\xi_i$ yields a map $\xi: K\cap i(G) \rightarrow \Z^k$ where the $j$-th coordinate is given by $\xi_j$.  We then obtain the following short exact sequence
$$1 \rightarrow Ker(\xi) \rightarrow K\cap i(G) \stackrel{\xi}\rightarrow \Z^k \rightarrow 1.$$
We would like to show that given any $\ph \in Aut(G)$, the kernel $Ker(\xi)$ is invariant under $\eta'_{\ph}$.  In terms of the quasi-isometry $\Theta'(g)|_{T_i'}$, we see that the image of the set of points at any height $h$ in $T_i'$ maps to a set of points of uniformly bounded height.

%
%

The kernel of the map $\xi$ is the set of $g \in K\cap i(G)$ so that the induced quasi-isometry $\Theta'(g_i)|_T$  preserves heights on each tree factor up to an additive constant. 
Alternately, if $g \in Ker(\xi)$ then the orbits $\{\Theta'(g_i)|_{T_i'}^n \cdot x\}$ stay at uniformly bounded height for each $i$.

Any  $\zeta \in Aut(K\cap i(G))$ is a quasi-isometry of $K\cap i(G)$ and thus induces a quasi-isometry on the product of trees $\ptp$.  A quasi-isometry will take sets of uniformly bounded height to sets of uniformly bounded height, and thus $\zeta(g)$ will again be in the kernel of $\xi$, and we see that $Ker(\xi)$ is invariant under $\eta'_{\ph}$ for any $\ph \in Aut(G)$.  Denote the restriction of $\zeta$ to $Ker(\xi)$ by $\zeta'$.  Thus we have obtained the vertical maps on the short exact sequence above.
\begin{equation}\label{general-exact}
\begin{CD}
    1 @>>> Ker(\xi)   @>>>  K\cap i(G) @>>>    \Z^k @>>> 1 \\
    @.     @V{\zeta'}VV  @V{\zeta}VV   @V{\overline \zeta}VV @.\\
    1 @>>> Ker(\xi)       @>>>  K\cap i(G) @>>>    \Z^k @>>> 1 
 \end{CD}
\end{equation}

The map ${\overline \zeta}: \Z^k \rightarrow \Z^k$ is factor preserving by construction, and thus for each factor of $\Z$, the map ${\overline \zeta}|_{\Z} = \pm Id$.  However, since each $\psi_i$ is a quasi-isometry which preserves the point at infinity, we conclude that ${\overline \zeta}|_{\Z} = Id$.  Then $R({\overline \zeta}) = \infty$, and it follows from Lemma \ref{lemma:reid} that $R(\zeta) = \infty$ as well.
\end{proof}

We now return to the short exact sequence 
\begin{equation}\label{general-exact}
\begin{CD}
    1 @>>> K\cap i(G)   @>>>  i(G) @>>>    F @>>> 1 \\
    @.     @V{\eta'_{\ph}}VV  @V{\eta_{\ph}}VV   @VVV @.\\
    1 @>>> K\cap i(G)   @>>>  i(G) @>>>    F @>>> 1 
 \end{CD}
\end{equation}

The map $\zeta$ in the proof of Lemma \ref{lemma:phi'} is the map $\eta'_{\ph}$ in this commutative diagram.  We have shown above that $R(\eta'_{\ph}) = \infty$, and it follows from Lemma \ref{lemma:reid} that $R( \eta_{\ph}) = \infty $ as well. Finally, we apply the same argument to $\overline{\ph}=\eta_{\ph}$ in diagram \eqref{center-exact} and conclude that $R(\varphi)=\infty$. 
\end{proof}


\begin{thebibliography}{FM2}

\bibitem[BMMV]{mv}
O. Bogopolski, A. Martino, O. Maslakova, and E. Ventura. The conjugacy problem for free-by-cyclic groups.  {\it Bull. London Math. Soc.}, to appear.

\bibitem[dlH]{dlH}
P. de la Harpe, \underline{Topics in Geometric Group Theory}, The University of Chicago Press, Chicago, IL (2000).

\bibitem[FM1]{fm1}
B. Farb and L. Mosher (appendix by D. Cooper), A rigidity theorem for
the solvable Baumslag-Solitar groups, {\it Inventiones Math},
Vol. 131, No. 2 (1998), 419-451.

\bibitem[FM2]{fm2}
B. Farb and L. Mosher , Quasi-isometric rigidity for the solvable Baumslag-Solitar groups, II, {\it Inventiones Math.}, Vol. 137, No. 3 (1999), 273-296.

\bibitem[F]{f}
A.L. Fel'shtyn, The Reidemeister number of any automorphism of a Gromov hyperbolic group is infinite, {\it Zap. Nauchn. Sem. POMI} Vol. 279 (2001), 229-241.

\bibitem[FH]{fh}
A.L. Fel'shtyn and R. Hill, The Reideemeister zeta function with applications to Nielsen theory and a connection with Reidemeister torsion, {\it K-theory}, Vol. 8 (1994), 367-393.

\bibitem[FG]{fg}
A.L. Fel'shtyn and D. Gon\c calves, Twisted conjugacy classes of automorphisms of Baumslag-Solitar groups, {\it Algebra and Discrete Math.}, to appear.

\bibitem[GW1]{gw1}
D. Gon\c calves and P. Wong, Twisted conjugacy in exponential growth groups, {\it Bull. London Math. Soc.}, Vol. 35 (2003), 261-268.

\bibitem[GW2]{gw2}
D. Gon\c calves and P. Wong, Twisted conjugacy classes in wreath products, {\it Internat. J. Alg. Comput.}, to appear.

%

\bibitem[K]{k}
D.H. Kouchloukova, The automorphism group of the Baumslag-Solitar group, preprint 2004. 

\bibitem[LL]{ll}
G. Levitt and M. Lustig, Most automorphisms of a hyperbolic group have simple dynamics, {\it Ann. Sci. Ecole Norm. Sup.} Vol. 33 (2000), 507-517.

\bibitem[M]{m}
J. Milnor, A note on curvature and fundamental group, {\it J. Differential Geom.} Vol. 2 (1968), 1-7.

\bibitem[MSW]{msw}
L. Mosher, M. Sageev and K. Whyte, Quasi-actions on trees I, {\it Ann. of Math.}(2) Vol. 158 (2003), 115-164.

\bibitem[S]{s}
Z. Sela,  Structure and rigidity in (Gromov) hyperbolic groups and discrete groups in rank $1$ Lie groups, II, {\it Geom. Funct. Anal.} Vol. 7, No. 3 (1997), 561-593. 


\bibitem[SW]{sw}
E. Souche and B. Wiest, An elementary approach to quasi-isometries of tree $\times \R^n$, {\it Geom. Dedicata} Vol. 95 (2002), 87-102.

\bibitem[TW]{tw}
J. Taback and K. Whyte, The Large Scale Geometry of Some Metabelian Groups, {\it Michigan  Math. Journal}, Vol. 52, No. 1 (2004), pp. 205-218.

\bibitem[TWo]{two}
J. Taback and P. Wong, A note on Twisted conjugacy and generalized Baumslag-Solitar groups, preprint, 2006.

\bibitem[Wh]{wh}
K. Whyte, Large scale dynamics, coarse fibrations, and the geometry of finitely generated groups, preprint, 2003.

\bibitem[W]{w}
P. Wong. Reidemeister number, Hirsch rank, coincidences on polycyclic groups and solvmanifolds. {\it J. reine angew. Math.}, Vol. 524 (2000), 185-204.
\end{thebibliography}
\end{document}